\documentclass{article}
\title{A sufficient condition for global invertibility of Lipschitz
mapping} \author {Sergey P.  Tarasov \thanks{Supported in part by RFBR
grant 02-01-00716} \\ Computing center of RAS \\ Vavilova 40 \\ 117967
Moscow GSP-1, Russia.\\ {\tt e-mail: sergey@ccas.ru} } \date{ }
  \newtheorem{theorem}{Theorem}

  \newtheorem{rK}{Remark}

  \newtheorem {pF}{Proof}

\let\ra\rightarrow
\let\eps\varepsilon
\def\RR{{\bf R}}

\def\J{{\cal J}}
\def\K{{\cal K}}
\def\M{{\cal M}}

\def\P{{\cal P}}

\begin{document}
\maketitle

\begin{abstract}
We show that S.Vavasis' sufficient
condition for
global invertibility of a polynomial mapping can be easily generalized
to the case of a general Lipschitz mapping.

\end{abstract}

In applications, in particular in finite-element analysis it is useful
to have some sufficient conditions for testing invertibility or
injectivity of a mapping $f : \RR^n \ra \RR^n$.
Usually $f(\cdot)$ is defined on some simple region $I$, say cube or
simplex, in the reference domain. The image $f(I)$ is some region with
more complicated geometry (a grid cell) in the physical domain.

In a paper \cite{Vavasis}
S.Vavasis proposed a sufficient
condition for global invertibility of a polynomial mapping in $\RR^n$.
Here we show that this condition can be more naturally formulated in
the  framework of {\em nonsmooth analysis} (see, e.g. \cite{Clarke}),
and this enables to spread the results of \cite{Vavasis} and
translate them almost word-for-word into a seemingly
much more general setting (at least in the linguistic sense).  Thus
this note may be regarded as a feedback to \cite{Vavasis}.

We need several definitions. The standard
reference on nonsmooth analysis is \cite{Clarke}.
Our discussion is restricted to a
finite-dimensional case though an essential part of \cite{Clarke} is
devoted to the infinite-dimensional generalizations.

Recall that nonsmooth analysis works with Lipschitz functions that
are almost everywhere differentiable.
Roughly speaking, the {\em generalized (or Clarke's)} gradient
$\partial F$ of a
Lipschitz function $f(\cdot): R^n \rightarrow R$ at the point $x \in R^n$
is defined  as a convex hull of (almost) all converging sequences of the
gradients: $\partial f(x) \stackrel{def}{=} conv(\lim \nabla(f(x_i))$,
for $x_i \rightarrow x$ and $f(\cdot)$ is differentiable at points $x_i$
and the sequence $\nabla f(x_i)$ converges. It is essential that
at the points of smoothness of $f(\cdot)$ the generalized gradient
coincides with gradient, and for a convex function--- with its subgradient.

Similarly
\footnote{Actually, Clarke's definitions of the generalized gradient
and the generalized Jacobian are somewhat less restrictive. We may
assume that the points of the sequences ${x_i}$ or $JF(x_i)$ involved
do not belong not only to nonsmooth point of the map but additionally
{\em they do not belong to an arbitrary set of Lebesque measure zero}.
Such definitions are technically more convenient.}, the {\em
generalized Jacobian} $\partial F$ of a Lipschitz mapping $F(\cdot):
R^n \rightarrow R^m$ at the point $x$ is a convex hull of all $m \times n$
matrices obtained as limits of sequences $JF(x_i)$
(its Jacobian matrix at the point of smoothness $x_i$),
where $x_i \rightarrow x$ and $F(\cdot)$ is differentiable at
$x_i$.

The generalized Jacobian {\em has maximum rank at point $x_0$}
if each matrix from $\partial F(x_0)$ has maximum rank.

Nonsmooth inverse function theorem \cite[Th.7.1.1]{Clarke}
states that Lipschitz mapping $F: R^n \rightarrow R^n$, whose
generalized Jacobian has maximum rank at point $x_0$ is locally
Lipschitz invertible.

In \cite{Vavasis} S.Vavasis proposed a sufficient
condition for the nondegeneracy of a  matrix family $\M$
consisting of
square $n \times n$ matrices. To present this
result in a slightly more general form we need some (more or
less standard) definitions.

{\em A cone} $\K \subseteq \RR^n$ is a set with a property that for all
$x \in \K$ and for any $\lambda \geq 0, \; \lambda x \in \K$, i.e. $\K$
contains all intersecting rays through the origin.

A cone $\K$ is {\em convex} if sum of any points of $\K$ belongs to
$\K$.

A convex cone is {\em polyhedral} it can be represented in the form:
$\K=\{x \,|\, x= \sum_{i=1}^m{\lambda_i a_i}, \; \lambda_i \geq 0, \;
a_i \in \RR^n \}$.

A convex cone $\K$ is {\em acute} if it contains no nontrivial
subspaces or equivalently, if no finite set of elements of $\K$ sums to
zero.

If $A \subseteq \RR^N$ then $cone\; hull(A)$ is a union of all rays
through the origin intersecting $A$.

Let formulate {\em Vavasis' condition}.
Denote by $R_i
\subseteq \RR^n$
the set of all $i$-th columns of the matrices
belonging to $\M$ and assume that all $R_i, \; i=1,\dots, n$ are
separated from the origin, i.e. for some $\delta>0 \; R_i \cap \{x \,
| \, \|x\| \leq \delta\} = \emptyset$ (here $\| \cdot \|$ denotes the
Euclidean norm).  Set $K_i = cone \; hull(R_i),\; i=1,\dots,n$ and let
any cone out of the $2^{n-1}$ cones $K_1 \pm K_2 \pm \dots \pm
K_n$ be {\em acute}.

By definition, $V$-{\em family} is any matrix family satisfying these
conditions. It follows from the above that all matrices in $V$-family
are nondegenerate.

If {\em additionally}, any cone involved
$K_1 \pm K_2 \pm \dots \pm K_n$
is contained in a cone $K_a=\{x
\in \RR^n \,| \, ax \geq \varepsilon \| x \|\}$ for some {\em
certificate vector} $a\in \RR^n$ and $\eps > 0$
then such family is defined as {\em strict
$V$-family}. In particular, if all $K_i$ (or $R_i$), $i=1,\dots, n$ are
polyhedral then by Farkas lemma all $V$-families are strict. Moreover,
checking that some matrix family is $V$-family is reduced to solving
$2^{n-1}$ linear programs, and thus the overall test could be performed
in linear time with respect to the input in any fixed dimension.

\begin{theorem} Let $F(\cdot)=(f^1(\cdot),\dots,f^n(\cdot)): U
\subseteq \RR^n \ra \RR^n$ be any Lipschitz mapping defined on the
convex reference domain $U$. If the set of the generalized Jacobians
$\J=\{J \in \partial F(u),\, u \in U\}$ forms a strict $V$-family then
$F(\cdot)$\ is globally invertible on $U$.  \end{theorem}

{\bf Proof.} Actually the demonstration is a direct
translation into the new setting of the original proof from
\cite{Vavasis}.

At first, local invertibility of $F(\cdot)$ at any point $u \in U$ of
the reference domain follows from the nonsmooth inverse function
theorem as the generalized Jacobian $\partial F(u)$ has
maximum rank at any point by construction.

Secondly, to show global invertibility it is enough to check
injectivity of the mapping. Take any different points $u, v  \in U$
from the reference domain. By convexity of $U$ and by the nonsmooth
analog of the Lagrange formula \cite[Th.2.6.5]{Clarke} for almost
all pairs $u, v \in U$ the following equality holds:
$F(v)-F(u)=\int^1_0{JF(u + t(v-u))(v-u)dt}$.
Assume w.l.o.g.  that the first
coordinate of the vector $v-u$ is nonnegative (otherwise, exchange $v$
and $u$). Set $K_i= cone\; hull(\partial f^i(x), \, x \in U)$. Now
assume for a moment that $v-u \geq 0$.  Take the corresponding
certificate vector $a$ for the cone $K_1+\dots+K_n$. By
construction, for all $s \in U, \; (a, \partial f_i(s)) > \eps \delta,
\; i=1,\dots,n$. Hence, $(a,F(v)-F(u)) > \eps \delta
\|u-v\| > 0$ $ $ and the injectivity follows. In general case, we
take any certificate vector for the cone $K_1 \pm K_2 \pm \dots
\pm K_n$ {\em with the same sign pattern as the sign pattern of the
coordinates of the vector $v-u$}.
\medskip

Informally, Vavasis sufficiency condition assumes that the columns of
the corresponding Jacobian matrices are independent and thus may be
extremely restrictive
but, on the other hand, it has attractive enough {\em decomposition
feature}.  Namely, let the reference domain $I$ be subdivided into
several simple regions, say, $I$ is a cube that is partitioned into
parallelepipedal patches $I=\cup_{i=1}^M I^{i}$ by several axis-
parallel hyperplanes. Now let $F(\cdot)=(f^1(\cdot),\dots,f^n(\cdot)):
I \ra \RR^n$ be any Lipschitz mapping.  Let $\Delta_i= \cup_{x \in
I_i}JF(x)$ (the union of all generalized Jacobians in the {\em closure}
of the patch $I_i$).  Obviously these sets can be computed {\em
separately} for any {\em closed} patch $I_i \; i=1, \dots,
M$. As above we obtain the following sufficient condition: if each of
the sets $\Delta_i, \, i=1, \dots, M$ forms a strict $V$-family and
there exists a certificate vector {\em common to all} $\Delta_i$ then
the mapping $F(\cdot)$ is globally invertible.  Equivalently, as
Vavasis condition is insensible to taking convex hulls, $F(\cdot)$ is
globally invertible on $I$ if the set $\Delta=conv(\cup_{i=1}^M
\Delta_i)$ forms a strict $V$-family.

For a simple example, assume
that $F(\cdot)$ is a continuous piecewise polynomial (product)
Bernstein-Bezier (BB) mapping as proposed in
\cite{Vavasis}\footnote{The reader should take into account that the
example below is almost explicit in \cite{Vavasis} and could be easily
recovered from the arguments therein {\em but with some additional
formal arguments.}} , i.e.  in any patch $I_q, \, q=1,\dots,
M$ the mapping is given by
$$F(\xi_1,\dots,\xi_n)=\sum_{i_1=0}^p\dots\sum_{i_n=0}^p f^q_{i_1,\dots
i_n} {p\choose i_1}\xi_1^{i_1}(1-\xi_1)^{p-i_1} \dots {p\choose
i_1}\xi_n^{i_n}(1-\xi_n)^{p-i_n},$$
where $(\xi_1,\dots,\xi_n)$ are
local coordinates in $I_q=\{0\leq \xi_1 \leq 1, \dots, 0\leq \xi_n \leq
1\}$ and the vectors $f^q_{i_1,\dots i_n} \in \RR^n, \,
i_1,\dots,i_n=1,\dots,n$ are called the {\em control points}. It is
well known that the image $F(I_q)$ of BB mapping is contained in the
convex hull of the control points. It is also well known that the
partial derivatives of $F(\cdot)$ can be put into BB form (and the
resulting control points for the derivatives are effectively computable
linear combinations of the control points for $F(\cdot)$). Thus for any
$i=1,\dots,n$ the set of all $i$-th columns of the Jacobian matrices
$JF(I_q)=\{JF(x), \; x \in I_q\}$ is contained in some {\em explicitly
computable polytope} $P^q_i$. Hence, all possible $i$-th columns of the
Jacobian matrices $JF(x),\; x \in I$ are contained in
$P_i=conv\{P^1_i,\dots, P^M_i\}, \; i=1,\dots,n$. (Here we use the
continuity of $F(\cdot)$ as at any point belonging to the intersection
of some patches $x \in I_{i_1} \cap \dots \cap I_{i_k}$ the generalized
Jacobian $JF(x)$ is by construction contained in
$conv\{JF(I_{i_1},\dots,JF(I_{i_k})\}$.)  Thus if the matrix family
$\P$, whose $i$-th column belongs to $P_i, \; i=1,\dots,n$, form a
strict $V$-family then $F(\cdot)$ is globally invertible on the whole
domain $I$.

Here we show that invertibility test in \cite{Vavasis} is valid for a
much more ample class of mappings and thus should be
{\em simultaneously robust and restrictive} enough.  Therefore  it
would be nice if some new arguments would be applied to the
invertibility problem even in the simplest case of bilinear polynomials.

\end{document}